\documentclass[12pt]{amsart}

\newcommand{\LL}{{\mathbb L}}
\newcommand{\V}{{\mathcal V}}
\newcommand{\AV}{\operatorname{AV}}
\newcommand{\SB}{\operatorname{SB}}

\newcommand{\C}{{\mathbb C}}
\newcommand{\Q}{{\mathbb Q}}
\newcommand{\Qbar}{{\overline{\Q}}}
\newcommand{\kbar}{{\overline{k}}}

\newcommand{\Fbar}{{\overline{\F}}}
\newcommand{\Z}{{\mathbb Z}}

\newcommand{\Aff}{{\mathbb A}}
\newcommand{\PP}{{\mathbb P}}
\newcommand{\OO}{{\mathcal O}}

\newcommand{\F}{{\mathbb F}}

\newcommand{\End}{\operatorname{End}}
\newcommand{\Lie}{\operatorname{Lie}}

\newcommand{\Gal}{\operatorname{Gal}}

\newcommand{\Pic}{\operatorname{Pic}}

\newcommand{\GL}{{\operatorname{GL}}}

\newcommand{\isom}{\simeq}

\newcommand{\tensor}{\otimes}

\newcommand{\directsum}{\oplus}

\newtheorem{theorem}{Theorem}
\newtheorem{lemma}[theorem]{Lemma}

\theoremstyle{definition}

\theoremstyle{remark}
\newtheorem{rem}{Remark$\!\!$}

\begin{document}

\title[Grothendieck ring of varieties]
{The Grothendieck ring of varieties is not a domain}
\author{Bjorn Poonen}
\thanks{This research was supported by NSF grant DMS-9801104,
          and a Packard Fellowship.}
\address{Department of Mathematics, University of California, Berkeley, CA 94720-3840, USA}
\email{poonen@math.berkeley.edu}
\date{April 15, 2002}


\maketitle

\section{The Grothendieck ring of varieties}
\label{statement}

Let $k$ be a field.
By a $k$-variety we mean a geometrically reduced, 
separated scheme of finite type over $k$.
Let $\V_k$ denote the category of $k$-varieties.
Let $K_0(\V_k)$ denote the free abelian group generated 
by the isomorphism classes of $k$-varieties,
modulo all relations of the form $[X-Y] = [X] - [Y]$
where $Y$ is a closed $k$-subvariety of a $k$-variety $X$.
Here, and from now on, $[X]$ denotes the class of $X$ in $K_0(\V_k)$.
The operation $[X]\cdot [Y] := [X \times_k Y]$
is well-defined, and makes $K_0(\V_k)$ a commutative ring with $1$.
It is known as the {\em Grothendieck ring of $k$-varieties}.
A completed localization of $K_0(\V_k)$ is needed for the theory
of {\em motivic integration}, 
which has many applications: see~\cite{looijenga2000} for a survey.

Our main result is the following.

\begin{theorem}
\label{main}
Suppose that $k$ is a field of characteristic zero.
Then $K_0(\V_k)$ is not a domain.
\end{theorem}

\begin{rem}
We conjecture that the result holds also for fields $k$
of characteristic $p$.
But we use a result whose proof relies on resolution of singularities
and weak factorization of birational maps,
which are known only in characteristic zero.
\end{rem}

\section{Abelian varieties of $\GL_2$-type}
\label{GL_2-type}

If $A$ is an abelian variety over a field $k_0$,
and $k$ is a field extension of $k_0$,
then $\End_k(A)$ denotes the endomorphism ring of the base extension
$A_k:=A \times_{k_0} k$, that is, the ring of endomorphisms
defined over $k$.

\begin{lemma}
\label{abelianvariety}
Let $k$ be a field of characteristic zero,
and let $\kbar$ denote an algebraic closure.
There exists an abelian variety $A$ over $k$
such that $\End_k(A) = \End_\kbar(A) \isom \OO$,
where $\OO$ is the ring of integers of a number field 
of class number $2$.
\end{lemma}

Let us precede the proof of Lemma~\ref{abelianvariety}
with a few paragraphs of motivation.
Our strategy will be to find a single abelian variety $A$ over $\Q$
such that the base extension $A_k$ works over $k$.

Suppose that $A$ is a nonzero abelian variety over $\Q$.
Let $\Lie A$ be its Lie algebra,
which is a $\Q$-vector space of dimension $\dim A$.
If $\End_\Q(A)$ is an order in a number field $F$,
then the $\OO$-action makes $\Lie A$ a vector space 
over $\OO \tensor \Q = F$; hence $[F:\Q] \le \dim_\Q \Lie A = \dim A$.
If moreover equality holds, then $A$ is said to be of {\em $\GL_2$-type}.
(The terminology is due to the following:
If $A$ is of $\GL_2$-type,
then the action of the Galois group $\Gal(\Qbar/\Q)$
on a Tate module $T_\ell A$ can be viewed as a representation
$\rho_\ell: \Gal(\Qbar/\Q) \rightarrow \GL_2(\OO \tensor_\Z \Z_\ell)$.)

Because $\Q$ has class number $1$,
we must take $[F:\Q] \ge 2$ 
to find an $A$ over $\Q$ as in Lemma~\ref{abelianvariety}.
The inequality $[F:\Q] \le \dim A$
then forces $\dim A \ge 2$.
Moreover, if we want $\dim A=2$,
then $A$ must be of $\GL_2$-type.

Abelian varieties of $\GL_2$-type are closely connected to modular forms.
For each $N \ge 1$, let $\Gamma_1(N)$ denote the classical modular group,
let $X_1(N)$ denote the corresponding modular curve over $\Q$,
and let $J_1(N)$ be the Jacobian of $X_1(N)$.
G.~Shimura \cite[Theorem~7.14]{shimura1971}, \cite{shimura1973} 
attached to each weight-2 newform $f$ on $\Gamma_1(N)$
an abelian variety quotient $A_f$ of $J_1(N)$.
It is known that $\dim A_f = [F:\Q]$,
where $F$ is the number field generated over $\Q$ 
by the Fourier coefficients of $f$.
These coefficients can also be identified
with the endomorphisms of $A_f$ 
induced by the Hecke correspondences on $X_1(N)$;
hence $A_f$ is $\GL_2$-type.
Conversely,
it is conjectured that each abelian variety of $\GL_2$-type is
$\Q$-isogenous to some $A_f$.
See~\cite{ribet1992} for more details.
The $\dim A=1$ case of this conjecture is the 
statement that elliptic curves over $\Q$ are modular,
which is known~\cite{breuil2001}.

Therefore we are led to consider $A_f$ of dimension $2$,
where $f$ is a newform as above.

\begin{proof}[Proof of Lemma~\ref{abelianvariety}]
Tables~\cite{steintables} show that
there exists a weight-2 newform $f= \sum_{n=1}^{\infty} a_n q^n$
on $\Gamma_0(276)$ (hence also on $\Gamma_1(276)$)
such that $\Q(\{a_n : n \ge 1\}) = \Q(\sqrt{10})$, 
$a_{17}=4-\sqrt{10}$, and $a_{19}=2+\sqrt{10}$.
Let $A=A_f$ be the corresponding abelian variety over $\Q$.
Then $\dim A = [\Q(\sqrt{10}):\Q] = 2$.
Also, $\End_\Q(A)$ is an order of $\Q(\sqrt{10})$ 
containing $4-\sqrt{10}$,
so $\End_\Q(A)$ is the maximal order $\Z[\sqrt{10}]$ of $\Q(\sqrt{10})$.
The class number of $\Q(\sqrt{10})$ is $2$.

It remains to show that $\End_k(A)=\Z[\sqrt{10}]$
for any field extension $k$ of $\Q$.
For any place of $k$ at which $A$ has good reduction,
$\End_k(A)$ injects into the endomorphism ring of the reduction.
We will use this to bound $\End_k(A)$.
The abelian variety $A$ has good reduction at all primes
not dividing $276$, so in particular it has good reduction at $17$ and $19$.
Let $A_{17}$ and $A_{19}$ denote the resulting abelian varieties
over $\F_{17}$ and $\F_{19}$.
The places $17$ and $19$ of $\Q$
extend to places of $k$ taking values in $\Fbar_{17}$ and $\Fbar_{19}$.
Thus $\End_k(A)$ injects into $\End_{\Fbar_{17}}(A_{17})$
and $\End_{\Fbar_{19}}(A_{19})$.

By the work of Eichler and Shimura (see Theorem~4 in D. Rohrlich's article
in~\cite{cornell-silverman-stevens1997}),
the characteristic polynomial $P_{17}(x)$ of Frobenius on $A_{17}$
equals 
	$$N_{\Q(\sqrt{10})/\Q} (x^2 - a_{17} x + 17)
	= x^4 - 8 x^3 + 40 x^2 - 136 x + 289.$$
This is irreducible over $\Q$, and its middle coefficient is prime to $17$,
so $A_{17}$ is a simple ordinary abelian surface.
Checking the criterion in~\cite{howe-zhu2002} (see especially 
Theorem~6 and the last paragraph of Section~2),
we find that $\End_{\Fbar_{17}}(A_{17}) \tensor \Q \isom \Q[x]/(P_{17}(x))$.
Similarly, $\End_{\Fbar_{19}}(A_{19}) \tensor \Q \isom \Q[x]/(P_{19}(x))$.
The ratio of the discriminants of $P_{17}(x)$ and $P_{19}(x)$
is not a square in $\Q$,
so $\Q[x]/(P_{17}(x))$ and $\Q[x]/(P_{19}(x))$ 
are {\em distinct} number fields of degree~4.
But $\End_k(A) \tensor \Q$ embeds into both,
so $\dim_\Q(\End_k(A) \tensor \Q) \le 2$.
On the other hand, $\Z[\sqrt{10}] \subseteq \End_k(A)$,
so $\End_k(A)=\Z[\sqrt{10}]$.
\end{proof}

\begin{rem}
The case $k=\C$ of Lemma~\ref{abelianvariety} has an easy proof:
let $A$ be an elliptic curve over $\C$ with 
complex multiplication by $\Z[\sqrt{-5}]$.
\end{rem}

\section{Abelian varieties and projective modules}
\label{projectivemodules}

Let $A$ be an abelian variety over a field $k$,
and let $\OO = \End_k(A)$.
Given a finite-rank projective right $\OO$-module $M$,
we define an abelian variety $M \tensor_\OO A$ as follows:
choose a finite presentation 
$\OO^m \rightarrow \OO^n \rightarrow M \rightarrow 0$,
and let $M \tensor_\OO A$ be the cokernel of the 
homomorphism $A^m \rightarrow A^n$ defined by the matrix
that gives $\OO^m \rightarrow \OO^n$.
It is straightforward to check that this is independent
of the presentation,
and that $M \mapsto (M \tensor_\OO A)$
defines a fully faithful functor $T$
from the category of finite-rank projective right $\OO$-modules
to the category of abelian varieties over $k$.
(Essentially the same construction is discussed in
the appendix by J.-P.~Serre in~\cite{lauter2001}.)

\begin{lemma}
\label{isomorphicsquares}
Let $k$ be a field of characteristic zero.
There exist abelian varieties $A$ and $B$ over $k$
such that $A \times A \isom B \times B$
but $A_{\kbar} \not\isom B_{\kbar}$.
\end{lemma}

\begin{proof}
Let $A$ and $\OO$ be as in Lemma~\ref{abelianvariety}.
Let $I$ be a nonprincipal ideal of $\OO$.
Since $\OO$ is a Dedekind domain,
the isomorphism type of a direct sum of fractional ideals $I_1 \directsum \dots \directsum I_n$ is determined exactly by the nonnegative integer $n$
and the product of the classes of the $I_i$ in the class group $\Pic(\OO)$.
Since $\Pic(\OO) \isom \Z/2$, 
we have $\OO \directsum \OO \isom I \directsum I$ as $\OO$-modules.
Applying the functor $T$ yields $A \times A \isom B \times B$,
where $B:= I \tensor_\OO A$.
Since $\End_{\kbar}(A)$ also equals $\OO$,
we have $B_{\kbar} = I \tensor_\OO A_{\kbar}$.
Since $T$ for $\kbar$ is fully faithful, $A_\kbar \not\isom B_\kbar$.
\end{proof}

\section{Rings related to the Grothendieck ring of varieties}
\label{rings}

For any extension of fields $k \subseteq k'$, there is a ring homomorphism
$K_0(\V_k) \rightarrow K_0(\V_{k'})$ mapping $[X]$ to $[X_{k'}]$.

Let $k$ be a field of characteristic zero.
Smooth, projective, geometrically integral $k$-varieties $X$ and $Y$
are called {\em stably birational}
if $X \times \PP^m$ is birational to $Y \times \PP^n$
for some integers $m,n \ge 0$.
The set $\SB_k$ of equivalence classes of this relation
is a monoid under product of varieties over $k$.
Let $\Z[\SB_k]$ denote the corresponding monoid ring.

When $k=\C$, there is a unique ring homomorphism 
$K_0(\V_k) \rightarrow \Z[\SB_k]$
mapping the class of any smooth projective integral variety
to its stable birational class~\cite{larsen-lunts2001}.
(In fact, this homomorphism is surjective, and its kernel
is the ideal generated by $\LL:=[\Aff^1]$.)
The proof in~\cite{larsen-lunts2001} requires resolution of singularities
and weak factorization of birational 
maps~\cite[Theorem~0.1.1]{abramovich2000}, 
\cite[Conjecture~0.0.1]{wlodarczyk2001}.
The same proof works over any algebraically closed field
of characteristic zero.

The set $\AV_k$ of isomorphism classes of abelian varieties over $k$
is a monoid.
The Albanese functor mapping a smooth, projective, geometrically integral
variety to its Albanese variety
induces a homomorphism of monoids $\SB_k \rightarrow \AV_k$,
since the Albanese variety is a birational invariant,
since formation of the Albanese variety commutes with products,
and since the Albanese variety of $\PP^n$ is trivial.
Therefore we obtain a ring homomorphism $\Z[\SB_k] \rightarrow \Z[\AV_k]$.

\section{Zerodivisors}
\label{zerodivisors}

\begin{proof}[Proof of Theorem~\ref{main}]
Let $A$ and $B$ be as in Lemma~\ref{isomorphicsquares}.
Then $([A]+[B])([A]-[B])=0$ in $K_0(\V_k)$.
On the other hand, $[A]+[B]$ and $[A]-[B]$ are nonzero,
because their images under the composition
	$$K_0(\V_k) \rightarrow K_0(\V_{\kbar}) \rightarrow 
	\Z[\SB_{\kbar}] \rightarrow \Z[\AV_{\kbar}]$$
are nonzero.
(The Albanese variety of an abelian variety is itself.)
\end{proof}


\section*{Acknowledgements}

I thank Eduard Looijenga and Arthur Ogus for discussions.
The package {\tt GP-PARI} was used to perform the calculations in 
the last paragraph of the proof of Lemma~\ref{abelianvariety}.

\providecommand{\bysame}{\leavevmode\hbox to3em{\hrulefill}\thinspace}

\end{document}